# Performance Analysis of Sparse Recovery Models for Bad Data Detection and State Estimation in Electric Power Networks

Junwei Yang, Wenchuan Wu, *Senior Member, IEEE*, Weiye Zheng, Wenjun Xian, Boming Zhang, Fellow, *IEEE*

*Abstract*—This paper investigates the sparse recovery models for bad data detection and state estimation in power networks. Two sparse models, the sparse $L_1$-relaxation model (L1-R) and the multi-stage convex relaxation model (Capped-L1), are compared with the weighted least absolute value (WLAV) in the aspects of the bad data processing capacity and the computational efficiency. Numerical tests are conducted on power systems with linear and nonlinear measurements. Based on numerical tests, the paper evaluates the performance of these robust state estimation models. Furthermore, suggestion on how to select parameter of sparse recovery models is also given when they are used in electric power networks.

## I. INTRODUCTION

POWER system state estimation is one of the core functions of energy management system (EMS). The key role of state estimation is to provide the estimates of state variables as accurate as possible. The method of weighted least square (WLS) was first proposed by F.C.Schweppe [1, 2, 3]. After that, many other methods like fast decoupled static state estimation (FDSE) [4] were developed. However, one of the problems of WLS is that it has no ability to process bad data in the measurements set, which are caused by communication or metering failures or malfunctions. To solve this problem some robust estimators were proposed. Rousseeuw and Leroy [5] presented the method of least median of squares (LMS) and least trimmed squares (LTS) to suppress the influence of bad data by using measurements selectively. Miranda [6] proposed a state estimation method based on the concept of correntropy to denoise the bad data. Gastoni [7] presented a robust state estimator based on the maximum agreement among measurements, choosing a best estimation in many estimation states. Besides, weighted least absolute value (WLAV) [8] is also one of the most popular robust estimators with good performance in the aspect of robustness. The core idea for most of these methods is to minimize the value of a penalty function of residuals, which means each residual is considered as a whole and undividable.

However, when we take the bad data into account, the fact is that the residual can be further divided into two different parts: usual additive observation noise and abnormally large measurement errors caused by bad data. The main idea of dealing with the bad data for sparse recovery methods is divided the residuals into two parts as mentioned before; they actually distributed differently and large measurement errors usually form a sparse vector in the process of modeling. Recently, several sparse models were proposed and their superiority in robustness were theoretically proved using the mathematical methods of compress sensing (CS). [9] proposed a convex programming model with sparse vector constrained by $l_\infty$-norm to make error correction and prove the effectiveness of the model. [10] used the sparse L1-relaxation model (L1-R) to make the sparse vector recovery and solve the optimization problem in linear and nonlinear system. After these researches, [11] put forward a new multi-stage convex relaxation model (Capped-L1) to handle these kind of issues with sparse error vector. All of these sparse models are aimed to make the sparse recovery. However, when they are used in the power system, their effectiveness, parameter sensitivity and computational efficiency are not clear enough as far as we know.

In this paper, we will make a comparison for three robust state estimators: WLAV, L1-R and the Capped-L1 in the aspect of their bad data processing capacity, parameter sensitivity and computational efficiency. In addition, through the case studies, a suggestion on how to tune parameter for L1-R and Capped-L1 is also given. To our best knowledge, Capped-L1 is firstly introduced for power system state estimation in this paper, and the numerical tests show it has better ability to denoise bad data than the other two robust state estimation models. The remainder of the paper is organized as follows. In Section II, formulations of the three methods are reviewed briefly. Section III describes the case studies and compares the three models by numerical results. Section IV concludes the paper.

## II. ROBUST STATE ESTIMATION MODELS

In this section, we will introduce three robust state estimation models: WLAV, L1-R and Capped-L1.

In general, when there exists bad data in the measurements, the relationship between the measurements and the states variables can be represented as follow: [10, 12]

$$y = h(x) + v + e \qquad (1)$$

where $y$ denotes the raw measurement vector, $h(x)$ denotes a

This work was supported by the National Key Basic Research Program of China (Gant. 2013CB228203) and in part by the National Science Foundation of China (Gant. 51477083)

J.Yang, W.Wu, W. Zheng and B.Zhang are with the Department of Electrical Engineering, State Key Laboratory of Power Systems, Tsinghua University, Beijing 100084, China. W. Xian is with State Grid Qinghai Electric Power Company, Xining, China.

set of measurement functions, which can be linear or nonlinear, $x$ denotes the vector of state variables, $v$ denotes the vector of noise and $e$ denotes the vector of errors corresponding to bad data.

In traditional models (e.g. WLAV and WLS), we actually treat each residual as a whole and the ultimate goal is to minimize the function of the residuals, while in sparse models, each residual is divided into two parts which are represented in the formulation (1) as $v$ and $e$. The noise vector $v$ is assumed to be $N(0,\sigma)$ and independent in most cases; the error vector $e$ is a sparse vector, generally containing few non-zero values.

In mathematics, when a vector is sparse, we can use $l0$-norm to denote the number of non-zero values in the vector. Naturally, in an optimization problem to get a sparse vector, we easily come up with the idea to formulate the objective function as $l0$-norm minimization. However, $l0$-norm minimization was proved to be a NP problem which cannot be solved efficiently [13]. Hence, some relaxation functions of $l0$-norm should be used. The main difference between the L1-R model and Capped-L1 model is about the relaxation function.

Based on the descriptions above, three robust models are reviewed briefly as follows.

*Model 1: Weighted least absolute value (WLAV) model*

The weighted least absolute value model can be formulated as [8]

$$\min \quad J(x) = \sum_i |w_i(y_i - h(x_i))| \quad (2)$$

Where $i$ is the index of measurements, and $w_i$ is the weight associated with $i^{\text{th}}$ measurement.

For clarity, in this paper, we simply consider the unweighted situation and ignore the difference between measurements' precision. Thus, the formulation can be simply described as follow using the definition of $l1$-norm:

$$\min \quad J(x) = \|y - h(x)\|_1 \quad (3)$$

This model minimizes the $l1$-norm of the residual vector to get robust state estimations. According to the convex optimization theory [14], $l1$-norm minimization helps obtain sparse solutions, which means the residual vector we finally get in this model is fairly sparse.

*Model 2: Sparse L1-relaxation (L1-R) model*

The original sparse recovery model is

$$\begin{aligned} \min \quad & J(x,e) = \|y - h(x) - e\|_2 \\ s.t. \quad & \|e\|_0 \leq \varepsilon \end{aligned} \quad (4)$$

However, (4) was proved to be a NP problem [13].

Thereby, a Sparse L1-relaxation model can be formulated as [10]

$$\begin{aligned} \min \quad & J(x,e) = \|y - h(x) - e\|_2 \\ s.t. \quad & \|e\|_1 \leq \varepsilon \end{aligned} \quad (5)$$

where $\varepsilon$ is a positive small number.

In this model, $l1$-norm is used to approximate $l0$-norm to get a sparse error vector $e$.

For given $y$ and $\varepsilon$, according to Lagrange duality theory, the solution to (5) corresponds to the solution in the following problem:

$$\min \quad J(x,e) = \|y - h(x) - e\|_2 + \lambda \|e\|_1 \quad (6)$$

for some Lagrange dual variable $\lambda \geq 0$

The Lagrangian relaxation algorithm is used in practice since we are not able to get the exact value of $\varepsilon$ in the real power system. By solving (6), we can obtain state estimations and error estimations simultaneously. Case studies in section III show that $\lambda$ is a pivotal parameter affecting the performance of L1-R model.

*Model 3: Multi-stage convex relaxation (Capped-L1) model*

The Capped-L1 model can be formulated as: [11]

$$\begin{aligned} \min \quad & J(x^{(l)}, e^{(l)}) = \|y - h(x^{(l)}) - e^{(l)}\|_2 + \lambda \cdot \sum_i c_i^{(l)} \cdot |e_i^{(l)}| \\ & c_i^{(l)} = I(|e_i^{(l-1)}| \leq \alpha^{(l)}) \in [0,1] \end{aligned} \quad (7)$$

It is a multi-stage convex relaxation model and $c$ is a relaxation parameter. $\alpha^{(l)}$ is a threshold value changing with the decision variable $e^{(l-1)}$. In the process of solving the problem, the following iteration will be needed:

---

Initialize $c_i^{(1)} = 1$ for $i = 1, \ldots n$
For $l = 1, 2, \ldots$
Let
$$\hat{e}^{(l)} = \arg\min \|y - h(x^{(l)}) - e^{(l)}\|_2 + \lambda \cdot \sum_i c_i^{(l)} \cdot |e_i^{(l)}|$$
Let $c_i^{(l+1)} = I(|e_i^{(l)}| \leq \alpha^{(l)})$ ($i = 1, \ldots n$)
$\alpha^{(l+1)} = f(e^{(l)})$

---

Fig. 1 Multi-stage Convex Relaxation for L1 Sparse Regularization (Capped-L1)

In the iteration process, $f(e^{(l)})$ is a function of $e^{(l)}$ and the threshold $\alpha^{(l)}$ is changing in each iteration.

Note that the each iteration process in Capped-L1 model is exactly solving a L1-R problem, which may mean the accuracy of Capped-L1 is at least same to L1-R model. In this Capped-L1 model, a non-convex regularization condition is used instead of the convex $L_1$ regularization to approximate $L_0$ regularization. Tong Zhang [11] claimed this multi-stage convex relaxation was closer to $L_0$ regularization. We will test this model through the following case studies.

### III. CASE STUDIES FOR COMPARING DIFFERENT MODELS

#### A. State estimation using measurements from phasor measurement units (PMU)

PMU can accurately measure voltage vectors and current vectors in power networks. Based on the PMUs' measurements, the state estimation model is linear. Since of budget limitations, PMUs cannot make the whole system observable. Nonetheless, partial networks may be observable with only PMUs' measurement, e.g. all the substations above 220kV are installed with PMU in China. Therefore, as discussed in [16], we assume that there are enough voltage vector measurements and current vector measurements from PMUs in this case study.

In this setting, the measurement function $h(x)$ is linear. The relationship between measurements and the state variables for a power system with n buses can be formulated as follows:

$$\dot{U}_i = e_i + jf_i$$
$$\dot{I}_i = \left\{\sum_{j=1}^{N}(G_{ij}e_j - B_{ij}f_j)\right\} + j\left\{\sum_{j=1}^{N}(G_{ij}f_j - B_{ij}e_j)\right\} \quad (8)$$
$$\dot{I}_{ij} = [ge_i - (b+y_c)f_i - ge_j + bf_j] + j[(b+y_c)e_i + gf_i - be_j - gf_j]$$
$$\dot{I}_{ij}' = [\frac{g}{k^2}e_i - \frac{b}{k^2}f_i - \frac{g}{k}e_j + \frac{b}{k}f_j] + j[\frac{b}{k^2}e_i + \frac{g}{k^2}f_i - \frac{b}{k}e_j - \frac{g}{k}f_j]$$

where $e_i$ and $f_i$ are the real part and the image part of the voltage of node $i$, which form the state variable vector; $\dot{U}_i$ and $\dot{I}_i$ are the measurements of voltage and injection current of node $i$; $\dot{I}_{ij}$ is the branch current from bus $i$ to bus $j$ for line $ij$, and $\dot{I}_{ij}'$ represents branch current for transformer winding $ij$, $k$ is the turn ratio of transformer; $G_{ij}$ and $B_{ij}$ are the real part and image part of the element in i$^{th}$ row and j$^{th}$ column of the admittance matrix for a given power network; $g$ and $b$ are the conductance and susceptance of line $ij$; $N$ is the number of nodes connected to node $i$; $y_c = \frac{Y}{2}$ and Y is the admittance against ground of the line $ij$. Given a power system, the admittances and the ratio of transformers are known.

The performance of the three proposed models are compared through simulations on an IEEE 9-bus transmission system, an IEEE 30-bus transmission system an IEEE 57-bus transmission system, and an IEEE 118-bus transmission system [17]. The standard deviation of the voltage noise is 0.002 p.u. and the standard deviation of current noise is 0.001 p.u. The bad data rate is chosen to be 6% in the cases. And the threshold $\alpha^{(l)}$ is assigned to be $10^5 * \min(\hat{e}^{(l-1)})$ which is relatively reasonable chosen by the way of simulating. We measure the estimation error $\|\hat{x} - x_{true}\|_2$ as the criterion of precision. Each $\|\hat{x} - x_{true}\|_2$ value is an average of 100 trials in which the locations of bad data are randomized.

Furthermore, since the $l1$-norm function is a nonlinear function, when we solve the L1-R model and Capped-L1 model, we relax the absolute function as follows: [18]

$$L1-R \quad model:$$
$$\min \quad J(x,e) = \|y - h(x) - e\|_2 + \lambda(a+b)$$
$$s.t. \quad a+b-e=0; \quad a \geq 0, b \geq 0.$$
$$Capped-L1 \quad model: \quad (9)$$
$$\min \quad J(x,e) = \|y - h(x) - e\|_2 + \lambda \cdot \sum_i c_i \cdot |e_i|$$
$$s.t. \quad a_i \geq 0; \quad b_i \geq 0; \quad a_i + b_i - e_i = 0; \quad i = 1,2\ldots n.$$

Figure 2, figure 3 and figure 4 show the estimation errors against $\lambda$ when bad data rate is 0.06 for the IEEE 9-bus transmission system, IEEE 57-bus system and the IEEE 118-bus system. As a result of the independence of WLAV model from $\lambda$, the trend for WLAV line is horizontal. We can conclude that L1-R model with appropriate $\lambda$ is better than WLAV model in accuracy, especially in the cases of large systems. However, Capped-L1 model has little superiority over L1-R despite needing more iterations times. By analyzing the change trends of estimation errors with Lagrange dual variable $\lambda$, we can suggest that $\lambda$ locates in the range of 0.2~0.5 to get better performance for L1-R and Capped-L1 used in linear systems.

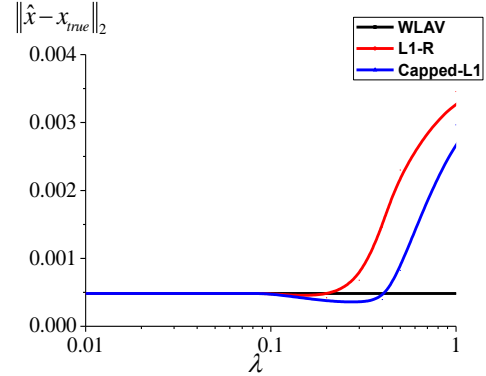

Fig. 2 Estimation error versus $\lambda$ for IEEE 9-bus system

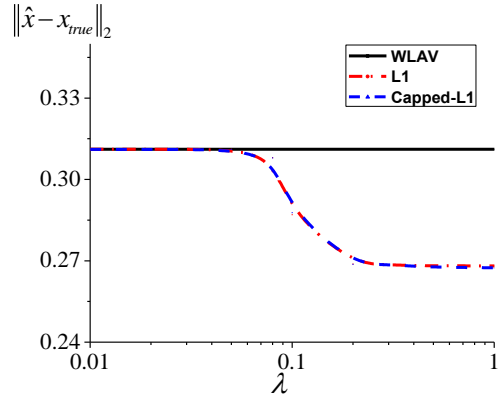

Fig. 3 Estimation error versus $\lambda$ for IEEE 57-bus system

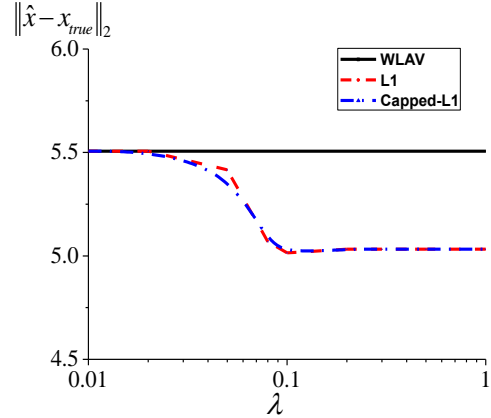

Fig. 4 Estimation error versus $\lambda$ for IEEE 118-bus system

We use the IBM Cplex to solve the optimization problem and Table1 shows the CPU time that three models take to estimate different system in the condition $\lambda = 0.3$. Obviously, WLAV is fastest and Capped-L1 is the slowest for IEEE 9-bus and 30-bus systems. However, for IEEE 57-bus and 118-bus, L1-R becomes the most efficient solution. Therefore, L1-R has advantages on both accuracy and efficiency with appropriate $\lambda$ for large systems.

TABLE 1
THE SOLVING TIME FOR THE ROBUST SE MODELS FOR LINEAR MEASUREMENTS

| Solving time/s | Case 9-bus | Case 30-bus | Case57-bus | Case118-bus |
|---|---|---|---|---|
| WLAV | 0.0054 | 0.0235 | 0.0576 | 0.2312 |
| L1-R | 0.0126 | 0.0268 | 0.0361 | 0.1094 |
| Capped-L1 | 0.0447 | 0.1217 | 0.1158 | 0.1789 |

Fig.5 shows the curve of estimation error against bad data percentage. All of the three models show good performance on denoising bad data with the bad data rate increasing. Capped-L1 and L1-R get similar estimation results with high precision, and they have superiorities over WLAV if an appropriate $\lambda$ is selected.

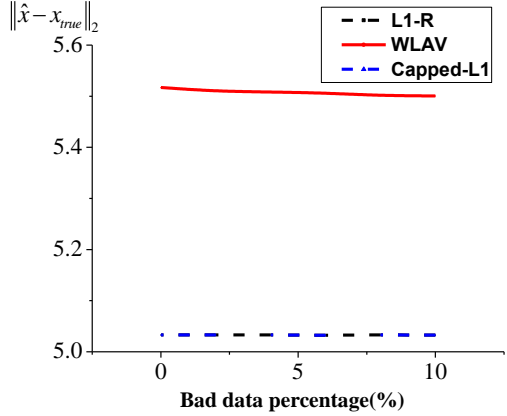

Fig. 5 Estimation error versus bad data percentage for IEEE 118-bus system with linear measurements

*B. State estimation using nonlinear measurements from supervisory control and data acquisition (SCADA) system*

In this part, we consider power system nonlinear state estimation with measurements from SCADA. The state variables are also the voltage magnitudes and the voltage angles of buses. The measurements include the injection real and reactive power at each bus and the real and reactive power flows on lines. In this case, the measurement function $h(x)$ is nonlinear and the optimization problem is a non-convex issue. A linearization process as shown in (10) is used to solve the state estimation models. [10]

linearized step of *WLAV* model:

$$\min \quad J(\Delta x) = \|\Delta y^k - H^k \Delta x\|_1$$

linearized step of $L1-R$ model:

$$\min \quad J(\Delta x, e) = \|\Delta y^k - H^k \Delta x - e\|_2 + \lambda(a+b)$$
$$s.t. \quad a \geq 0; \quad b \geq 0; \quad a+b-e = 0. \quad (10)$$

linearized step of $Capped-L1$ model:

$$\min \quad J(\Delta x, e) = \|\Delta y^k - H^k \Delta x - e\|_2 + \lambda \cdot \sum_i c_i \cdot |a_i + b_i|$$
$$c_i = I(|e_i| \leq \alpha) \in [0,1]$$
$$s.t. \quad a_i \geq 0; \quad b_i \geq 0; \quad a_i + b_i - e_i = 0.$$

where $\Delta y^k = y - h(x^k)$, $H^k$ is the Jacobian matrix of $h$ at $x^k$, $\Delta x$ is the solution to (10), and the state vector will be updated by

$$x^{k+1} = x^k + \Delta x \quad (11)$$

The iteration process is repeated until $|\Delta x| \leq 10^{-5}$, the final solution can be gotten.

In the nonlinear case studies, the standard deviation of the voltage noise is 0.002 p.u. and the standard deviation of power noise is 0.004 p.u. The threshold $\alpha^{(l)}$ is assigned to be $10^5 * \min(\hat{e}^{(l-1)})$. Similarly, we measure the estimation errors by $\|\hat{x} - x_{true}\|_2$.

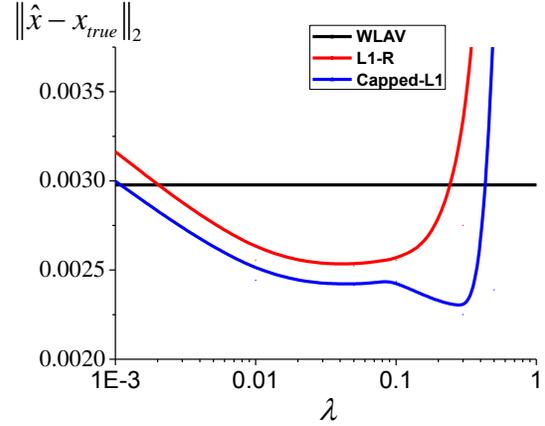

Fig. 6 Estimation error versus $\lambda$ for IEEE 9-bus system

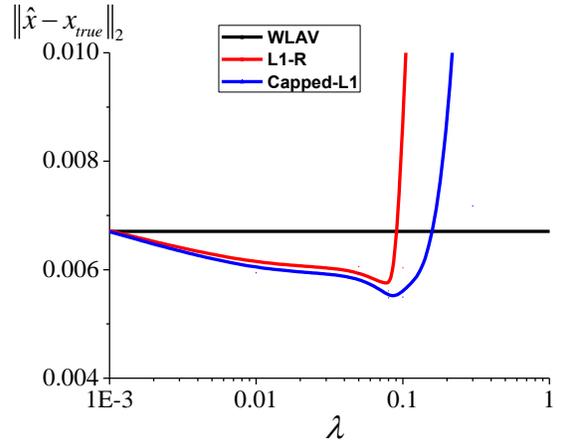

Fig. 7 Estimation error versus $\lambda$ for IEEE 30-bus system

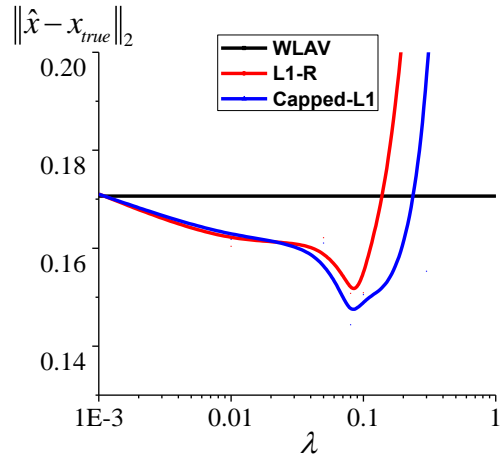

Fig. 8 Estimation error versus $\lambda$ for IEEE 57-bus system

Figure 4, figure 5 and figure 6 show an overall trend of error changes against $\lambda$ in the case that bad data rate is 0.06 for IEEE 9-bus, 30-bus and 57-bus system. L1-R model and Capped-L1 model is fairly sensitive to parameter $\lambda$, showing as a large error when $\lambda$ is large enough. In this condition, the superiority of L1-R model is not obvious, so we detail the con-

trast of estimation errors in the condition of $\lambda = 0.08$ in table 2. Table 2 lists the estimation errors and the percentages of the ascension on accuracy that two sparse recovery models versus WLAV model for each simulation cases. From the results listed in table 2, we can see L1-R model and Capped-L1 model are better than WLAV in processing bad data.

TABLE 2
THE STATE ERRORS FOR THE ROBUST SE MODELS FOR NONLINEAR MEASUREMENTS

| State Errors | Case 9-bus | Case 30-bus | Case57-bus | Case118-bus |
|---|---|---|---|---|
| WLAV | 0.0030 | 0.0067 | 0.1706 | 2.0076 |
| L1-R | 0.0026 | 0.0056 | 0.1508 | 1.9751 |
| Capped-L1 | 0.0024 | 0.0055 | 0.1444 | 1.9751 |
| Improvement ( L1-R versus WLAV) | 13.33% | 16.42% | 11.61% | 1.619% |
| Improvement (Capped-L1 versus WLAV) | 20.00% | 17.91% | 15.36% | 1.619% |

Table 3 lists the CPU times that three models take in each case when $\lambda = 0.08$. Similar to the conclusions made in the cases with linear state estimation, WLAV is fastest for small test system but L1-R becomes the most efficient for large system. Capped-L1 is the most time consuming as a result of additional iterations. Therefore, from table 2 and 3, we can claim L1-R model has advantages in both aspects of precision and computational efficiency if $\lambda$ locates in the range of 0.05~0.1. Whereas, Capped-L1 can get most accurate estimations.

TABLE 3
THE SOLVING TIME FOR THE ROBUST SE MODELS FOR NONLINEAR MEASUREMENTS

| Solving time/s | Case 9-bus | Case 30-bus | Case57-bus | Case118-bus |
|---|---|---|---|---|
| WLAV | 0.0059 | 0.0252 | 0.0721 | 0.2745 |
| L1-R | 0.0159 | 0.0655 | 0.0931 | 0.2384 |
| Capped-L1 | 0.0582 | 0.1774 | 0.2901 | 0.4241 |

Fig.9 shows that all the three models are robust enough for scenarios with different bad data percentages, in which $\lambda = 0.08$.

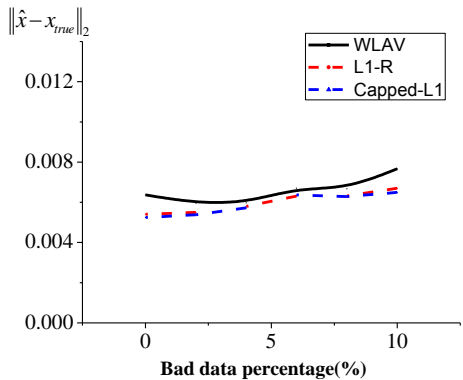

Fig. 9 Estimation error versus bad data percentage for IEEE 30-bus system with nonlinear measurements

## IV. CONCLUSIONS

In this paper, the performance of three robust state estimation models were compared and analyzed. Capped-L1 model is firstly introduced for power system robust state estimation here. Through case studies, we conclude that in the aspect of precision, Capped-L1 can get minimum estimation errors with appropriate parameter setting. However, Capped-L1 model has heavy computation burden. Estimation errors of L1-R model is slightly larger than those of Capped-L1 model, but L1-R has much higher computational efficiency. WLAV model is accurate and efficient but has no advantages compared to L1-R model. To summarize, L1-R model performs best in both denoising bad data and computation efficiency for large system, while Capped-L1 model has most powerful potential to process bad data with heavier computation burden. The Lagrange dual variable $\lambda$ in L1-R model and Capped-L1 is a vital parameter. By means of case studies, we suggest that $\lambda$ is assigned in the range of 0.2~0.5 for linear state estimation models and 0.05~0.1 for nonlinear state estimation models.